\documentclass[12pt]{article}
\usepackage{amsmath}
\usepackage{color}
\usepackage{fancyhdr}
\usepackage{verbatim}
\usepackage{amsfonts,longtable,mathrsfs}
\usepackage{epsfig}
\usepackage{amsfonts}
\usepackage{color}
\usepackage{fullpage}
\usepackage[numbers,sort&compress]{natbib}

\def\bn{\begin{equation}}
\def\en{\end{equation}}
\def\bny{\begin{eqnarray}}
\def\eny{\end{eqnarray}}
\def\be{\begin{eqnarray*}}
\def\ee{\end{eqnarray*}}
\def\bc{\begin{center}}
\def\ec{\end{center}}
\newtheorem{dfn}{Definition}[section]

\begin{document}

\bc {\bf On a system of second-order difference equations
  }\ec
\medskip
\bc
M. Folly-Gbetoula\footnote{Corresponding author: Mensah.Folly-Gbetoula@wits.ac.za}  and D. Nyirenda \footnote{Author: Darlison.Nyirenda@wits.ac.za} \vspace{1cm}
\\School of Mathematics, University of the Witwatersrand, Johannesburg 2050, South Africa.\\

\ec
\begin{abstract}
\noindent  We obtain explicit formulas for the solutions of the system of second-order difference equations of the form
\begin{align*}
x_{n+ 1}=\frac{x_{n}y_{n-1}}{y_{n}(a_n+b_nx_{n}y_{n - 1})}, \quad y_{n+1}=\frac{x_{n - 1}y_{n}}{x_{n}(c_n+d_nx_{n-1}y_{n})},
\end{align*}  where $(a_n)_{n\in \mathbb{N}_0},\;(b_n)_{n\in \mathbb{N}_0},\;(c_n)_{n\in \mathbb{N}_0}$ and $(d_n)_{n\in \mathbb{N}_0}$ are real sequences. We use Lie symmetry analysis to derive non-trivial symmetries  and thereafter, exact solutions are obtained.
\end{abstract}
\section{Introduction} \setcounter{equation}{0}
Lie symmetry analysis has been widely used to obtain solutions of systems of differential equations. This symmetry has diverse applications for instance the reduction of order of the equations via the invariants of their symmetry groups. Recently the method has been applied to difference equations and it is fruitful \cite{QR,hydon1,JV, W}. In applying the analysis to systems of difference equations, just as in the case of differential equations, one has to find a certain group of transformations that leave the equation invariant, but simply permute the set of solutions. Hydon \cite{hydon0} constructed a systematic methodology which can be used to find the group of transformations for difference equations. However, calculations can be cumbersome and to the best of our knowledge, there are no computer software packages that generate symmetries for difference equations. For ideas on how to derive solutions via the symmetry approach, the reader is referred to \cite{hydon1,FK3, FK,FK2}. \par \noindent
Our interest is in rational ordinary difference equations, which have been researched widely using different approaches, see \cite{ab-1,syst, ssys,RK,stevo,TE,IY}.
Elsayed, in \cite{syst}, solved the system:
\begin{align}\label{el}
x_{n+1}=\frac{x_{n-1}}{\pm y_nx_{n-1}\pm 1}, \quad y_{n+1}=\frac{ y_{n-1}}{\pm x_n
y_{n-1}\pm 1}, \quad n\in {N_0}.
\end{align}
In \cite{ssys}, the authors investigated \begin{align}\label{ei}
x_{n+1}=\frac{x_{n-2}y_{n-1}}{y_n( \pm 1 \pm y_{n-1}x_{n-2})}, \quad y_{n+1}=\frac{ y_{n-2}x_{n-1}}{x_n(\pm 1 \pm x_{n-1}
y_{n-2})}.
\end{align}
In this paper, we obtain explicit formulas for solutions of the system
\begin{align}\label{1.1}
x_{n+1}=\frac{x_{n}y_{n-1}}{y_{n}(a_n+b_nx_{n}y_{n-1})}, \quad y_{n+1}=\frac{x_{n-1}y_{n}}{x_{n}(c_n+d_nx_{n-1}y_{n})},
\end{align}  where $(a_n)_{n\in \mathbb{N}_0},\;(b_n)_{n\in \mathbb{N}_0}, \;(c_n)_{n\in \mathbb{N}_0}$ and $(d_n)_{n\in \mathbb{N}_0}$ are non-zero real sequences.
\subsection{Preliminaries}
A background on Lie symmetry analysis of difference equations is presented in this section. The notation used is that from \cite{hydon0}.
\begin{dfn}\cite{P.Olver} Let $G$ be a local group of transformations acting on a manifold $M$. A set $\mathcal{S}\subset M$ is called $G$-invariant, and $G$ is called symmetry group of $\mathcal{S}$, if whenever $x\in \mathcal{S} $, and $g\in G$ is such that $g\cdot x$ is defined, then $g\cdot x \in \mathcal{S}$.
\end{dfn}
\begin{dfn}\cite{P.Olver} Let $G$ be a connected group of transformations acting on a manifold $M$. A smooth real-valued function $\zeta: M\rightarrow \mathbb{R}$ is an invariant function for $G$ if and only if $$X(\zeta)=0\qquad \text { for all } \qquad  x\in M,$$
and every infinitesimal generator $X$ of $G$.
\end{dfn}
\begin{dfn}\cite{hydon0}
A parameterized set of point transformations,
\begin{equation}
\Gamma_{\varepsilon} :x\mapsto \hat{x}(x;\varepsilon),
\label{eq: b}
\end{equation}
where $x=x_i, $ $i=1,\dots,p$ are continuous variables, is a one-parameter local Lie group of transformations if the following conditions are satisfied:
\begin{enumerate}
\item $\Gamma_0$ is the identity map if $\hat{x}=x$ when $\varepsilon=0$
\item $\Gamma_a\Gamma_b=\Gamma_{a+b}$ for every $a$ and $b$ sufficiently close to 0
\item Each $\hat{x_i}$ can be represented as a Taylor series (in a neighborhood of $\varepsilon=0$ that is determined by $x$), and therefore
\end{enumerate}
\begin{equation}
\hat{x_i}(x:\varepsilon)=x_i+\varepsilon \xi _i(x)+O(\varepsilon ^2), i=1,\ldots,p.
\label{eq: c}
\end{equation}
\end{dfn}
Consider a given system of difference equations of order two:
\begin{align}\label{general}
\begin{cases}
x_{n+2}=&\Omega _1(n,x_n, x_{n+1}, x_{n+1}, y_n, y_{n+1}, y_{n+1}), \\
y_{n+k+1}=&\Omega _2(n,x_n, x_{n+1}, x_{n+1}, y_n, y_{n+1}, y_{n+1}), \quad n\in D
\end{cases}
\end{align}
for some smooth functions $\Omega_i$, $i=1,2$, and a domain $D\subset \mathbb{Z}$.
To compute a symmetry group of \eqref{general}, we pay attention to the group of point transformations given by
\begin{equation}\label{Gtransfo}
G_{\varepsilon}: (x_n,y_n) \mapsto(x_n+\varepsilon Q_1 (n,x_n),y_n+\varepsilon Q_2 (n,y_n)),
\end{equation}
where $\varepsilon$ is the parameter and $Q_i,\; i=1,2$, the continuous functions which we refer to as characteristics. Let
\begin{align}\label{Ngener}
X= & Q_1(n,x_n)\frac{\partial}{ \partial x_n}+ Q_2(n,y_n)\frac{\partial}{ \partial y_n}
\end{align}
be the  infinitesimal of $G _{\varepsilon}$. The group of transformations
$G_{\varepsilon}$ is a symmetry group of \eqref{general} if and only if
\begin{subequations}\label{LSC}
\begin{align}
\mathcal{S}^{(2)} Q_1- X^{[1]} \Omega _1=0\\
\mathcal{S}^{(2)} Q_2- X^{[1]} \Omega _2=0
\end{align}
\end{subequations}
where
\begin{equation}
X^{[1]}=
Q_1\frac{\partial}{ \partial x_n}+ Q_2\frac{\partial}{ \partial y_n}+ \mathcal{S}Q_1\frac{\partial\qquad }{ \partial x_{n+1}} +\mathcal{S}Q_2\frac{\partial\qquad}{ \partial y_{n+1}} + \mathcal{S}Q_1\frac{\partial\qquad}{ \partial x_{n+1}} +\mathcal{S}Q_2\frac{\partial\qquad}{ \partial y_{n+1}}
\end{equation}
because $\Omega _1$ and $\Omega _2$ are functions of $x_n, x_{n+1}, y_n$ and $y_{n+1}$ only.
The shift operator, $\mathcal{S}$, is defined as follows: $\mathcal{S}:n \rightarrow n+1$. Once we know the characteristics $Q_i$, the invariant $\zeta  _i$ can be deduced  by introducing the canonical coordinate \cite{JV}
 \begin{align}\label{cano}
s_n= \int{\frac{du_n}{Q_1(n,u_n)}} \quad \text{ and }\quad t_n= \int{\frac{du_n}{Q_2(n,u_n)}} .
 \end{align}
In general, the constraints on the constants in the characteristics  tell more the perfect choice of invariants, as opposed to lucky guesses.
\section{Symmetries and reductions}
Consider the system of difference equations
\begin{align}\label{1.2'}
\begin{cases}
x_{n+2}=\Omega _1 =\frac{x_{n+1}y_{n}}{y_{n+1}(A_n+B_nx_{n+1}y_{n})}\\ \\ y_{n+2}=\Omega _2 =\frac{x_ny_{n+1}}{x_{n+1}(C_n+D_nx_ny_{n+1})},
\end{cases}
\end{align}
 where $(A_n)_{n\in \mathbb{N}_0},\;(B_n)_{n\in \mathbb{N}_0}, \;(C_n)_{n\in \mathbb{N}_0}$ and $(D_n)_{n\in \mathbb{N}_0}$ are non-zero real sequences, equivalent to \eqref{1.1}.
\subsection{Symmetries}
To compute the symmetries, we impose condition \eqref{LSC} and obtain
\begin{subequations}\label{a1A1}
\begin{align}\label{a1}
& -Q_1(n+2,x_{n+2} )+ \frac{B_n {y_n}^{2}x_{n+1}^2Q_2(n+1,y_{n+1})}{{y_{n+1}}^2(A_n+B_nx_{n+1}y_{n})^2}-
\frac{A_n {y_n}y_{n+1}Q_1(n+1, x_{n+1})}{{y_{n+1}}^2(A_n+B_nx_{n+1}y_{n})^2}\nonumber\\
& +\frac{A_n {y_n}x_{n+1}Q_2(n+1,y_{n+1})}{{y_{n+1}}^2(A_n+B_nx_{n+1}y_{n})^2}-
\frac{A_nx_{n+1}y_{n+1}Q_2(n,y_n)}{{y_{n+1}}^2(A_n+B_nx_{n+1}y_{n})^2}=0,
\end{align}
\begin{align}\label{A1}
& -Q_2(n+2,y_{n+2})+ \frac{D_n {x_n}^{2}y_{n+1}^2Q_1(n+1, x_{n+1})}{(x_{n+1}(C_n+D_nx_ny_{n+1}))^2}-
\frac{C_n {x_n}x_{n+1}Q_2(n+1, y_{n+1})}{(x_{n+1}(C_n+D_nx_ny_{n+1}))^2}\nonumber\\
& +\frac{C_n {x_n}y_{n+1}Q_1(n+1, x_{n+1})}{(x_{n+1}(C_n+D_nx_ny_{n+1}))^2}-
\frac{C_ny_{n+1}x_{n+1}Q_1(n, x_{n})}{(x_{n+1}(C_n+D_nx_ny_{n+1}))^2}=0.
\end{align}
\end{subequations}
These functional equations for the characteristics $Q_i,\;i=1,2$ make \eqref{a1A1} hard to solve. We now eliminate the arguments $x_{n+2}$ and $y_{n+2}$ by operating the differential operators
\begin{subequations}
\begin{align}
L_1 =\frac{\partial}{\partial x_{n+1}}+\frac{\partial y_{n}}{\partial x_{n+1}}\frac{\partial}{\partial y_{n}}=\frac{\partial}{\partial x_{n+1}}-\frac{{\Omega_1} _{,x_{n+1}}}{{\Omega_1} _{,y_{n}}}\frac{\partial\quad}{\partial y_{n}}
\end{align}
on \eqref{a1} and
\begin{align}
L_2 =\frac{\partial}{\partial x_n}+\frac{\partial y_{n+1}}{\partial x_n}\frac{\partial}{\partial y_{n+1}}=\frac{\partial}{\partial x_n}-\frac{{\Omega_2} _{,x_n}}{{\Omega_2} _{,y_{n+1}}}\frac{\partial\quad}{\partial y_{n+1}}
\end{align}
\end{subequations}
on  \eqref{A1}. Note that $\Omega_{,x}$ is the partial derivative of $\Omega$ with respect to $x$. This results in
\begin{subequations}\label{a3A3}
\begin{align}\label{a3}
-y_{n}Q_1'(n+1,x_{n+1})+y_{n}Q_2'(n,y_n)+\frac{y_{n}}{x_{n+1}}Q_1(n+1,x_{n+1})-Q_2(n,y_n)=0
\end{align}
and
\begin{align}\label{A3}
y_{n+1}Q_2'(n+1,y_{n+1})-y_{n+1}Q_1'(n,x_n)+\frac{y_{n+1}}{x_n}Q_1(n,x_n)-Q_2(n+1,y_{n+1})=0
\end{align}
\end{subequations}
when fractions are cleared.
To eliminate the arguments $x_{n+1}$ and $y_{n+1}$, we divide both sides of \eqref{a3} by $y_n$ and differentiate with respect to $y_n$; differentiate \eqref{A3} with respect to $x_n$. Solving the resulting differential equations for $Q_1$ and $Q_2$ leads to
\begin{subequations}\label{a5A5}
\begin{align}\label{a5}
Q_1(n,x_n)=\alpha_n x_n + \beta _n x_n \ln x_n
\end{align}
and
\begin{align}\label{A5}
Q_2(n,y_{n})=\lambda _n y_n + \mu _n y_n\ln y_n,
\end{align}
\end{subequations}
where $\alpha _n,\; \beta _n,\; \lambda _n $ and $\mu _n$ depend of $n$ arbitrarily. The dependence among these functions is found by substituting equations in \eqref{a5A5} in equations in \eqref{a1A1}. The equations thereafter, can be solved by the method of separation which gives rise to the following systems:
\begin{subequations}\label{a8A8}
\begin{align}\label{a8}
\begin{cases}
x_{n+1}{y_n}&: \lambda _{n+1} +\alpha _{n+2}=0\\
1&: \lambda _{n+1} +\alpha _{n+2}-\alpha _{n+1} -\lambda _{n} =0
\end{cases}
\end{align}
and
\begin{align}\label{A8}
\begin{cases}
x_{n}y_{n+1}&: \alpha _{n+1} +\lambda _{n+2}=0\\
1&: \alpha _{n+1} +\lambda _{n+2}-\lambda _{n+1} -\alpha _{n} =0
\end{cases}
\end{align}
\end{subequations}
or simply
\begin{subequations}\label{a9A9}
\begin{align}\label{a9}
\lambda _n + \alpha _{n+1} =0
\end{align}
and
\begin{align}\label{A9}
\alpha _n +\lambda _{n+1} =0.
\end{align}
\end{subequations}
One can show that $\beta _n$ and $\mu _n$ are zero.
From \eqref{a9A9}, we note that
\begin{align}
\lambda_{n+2}-\lambda _n =0.
\end{align}\label{A9'}
Equation \eqref{A9'}  has
\begin{align}\label{A9''}
\lambda _n=c_0+(-1)^nc_1
\end{align}
 as general solutions
 and so, thanks to \eqref{a9A9} and \eqref{A9''}, the characteristics are;
\begin{align}\label{a10}
Q_{11}=x_n,\quad Q_{12}= (-1)^n x_n,\quad  Q_{21}=y_n,\quad Q_{22}=(-1)^ny_n.
\end{align}
Hence, the symmetry generators of \eqref{1.1} are
\begin{align}\label{a11}
X_1= & x_n\frac{\partial}{ \partial x_n}+ y_n\frac{\partial}{ \partial y_n}\qquad \qquad \quad
\end{align}
and
\begin{align}\label{a11'}
X_2= & (-1)^n x_n\frac{\partial}{ \partial x_n}+(-1)^n y_n\frac{\partial}{ \partial y_n}.
\end{align}
\subsection{Reductions}
Using \eqref{cano} and \eqref{a11'}, the canonical coordinates then found to be
\begin{align}\label{11}
s_n =(-1)^n\ln |x_n | \quad \text{and} \quad t_n=(-1)^{n}\ln| y_n|.
\end{align}
We replace $\alpha _n$ and its shift (resp $\lambda _n$ and its shift) with $s _n \alpha _n $ and its shift (resp $t_n \lambda _n $ and its shift) in \eqref{a9A9} and the left hand sides of the resulting equations give the invariants:
\begin{align}\label{a12}
\tilde{U}_n=\lambda _n t_n+ \alpha _{n+1}s_{n+1}=\ln|y_nx_{n+1}|
\end{align}
and
\begin{align}\label{A12}
\tilde{V}_n=\alpha _n s_n +\lambda _{n+1} t_{n+1}=\ln|x_{n}y_{n+1}|.
\end{align}
One can easily verify that $X[\tilde{U}_n]=X[\tilde{V}_n]=0$. For convenience, we use
\begin{align}\label{13}
U_n =\exp\{-\tilde{U}_n\}\quad \text{and} \quad
V_n=\exp\{-\tilde{V}_n\}
\end{align}
instead, or simply
\begin{align}\label{13'}
U_n=\pm\frac{1}{x_{n+1}y_{n}} \quad  \text{ and } \quad V_n=\pm\frac{1}{x_{n}y_{n+1}}.
\end{align}
Using the plus sign, this develops into
\begin{subequations}\label{14}
\begin{align}
    U_{n+1}  =  A_nU_n+B_n\label{14'}\\
    V_{n+1} =C_nV_n +D_n.\label{14''}
\end{align}
\end{subequations}
After iteration, it is easy to see that the solutions of equations in \eqref{14} in closed form are given by
\begin{subequations}\label{a16A16}
\begin{align}\label{a16}
U_{n}=&U_0 \prod\limits _{k_1=0}^{n-1}A_{k_1} +\sum\limits_{l=0}^{n-1} B_{l}\prod\limits_{k_2=l+1}^{n-1}A_{k_2},
\end{align}
\begin{align}\label{A16}
V_{n}=&V_0\prod _{k_1=0}^{n-1}C_{k_1} +\sum\limits_{l=0}^{n-1} D_{l}\prod\limits_{k_2=l+1}^{n-1}C_{k_2}.
\end{align}
\end{subequations}
Consequently, from \eqref{13'}, we have
\begin{align}\label{a17}
x_{n+2}=\frac{V_n}{U_{n+1}}x_n \quad \text{and} \quad  y_{n+2}=\frac{U_n}{V_{n+1}}y_n.
\end{align}
After some iterations, one finds that
\begin{align}\label{a18}
x_{2n+j}= x_j\prod_{s = 0}^{n - 1}\frac{V_{2s + j}}{U_{2s + j +1}} \quad \text{and} \quad  y_{2n+j}=y_j\prod_{s = 0}^{n - 1}\frac{U_{2s + j}}{V_{2s + j +1}}.
\end{align}
where $j = 0, 1$. So we have
\begin{align*}
x_{2n} & = x_0\prod_{s = 0}^{n - 1}\frac{    V_{2s}}{U_{2s +1}} \\
       & = x_0\prod_{s = 0}^{n - 1}\frac{V_0 \prod\limits _{k_1=0}^{2s-1}C_{k_1} +\sum\limits _{l=0}^{2s-1}  D_{l}\prod\limits _{k_2=l+1}^{2s-1}C_{k_2}}{ U_0 \prod\limits _{k_1=0}^{2s}A_{k_1} +\sum\limits _{l=0}^{2s}  B_{l}\prod\limits _{k_2=l+1}^{2s}A_{k_2}          } \\
       & = x_0 \frac{V_0^{n}}{U_0^{n}}\prod_{s = 0}^{n - 1}\frac{\prod\limits _{k_1=0}^{2s-1}C_{k_1} + \frac{1}{V_0}\sum\limits _{l=0}^{2s-1}  D_{l}\prod\limits _{k_2=l+1}^{2s-1}C_{k_2}}{\prod\limits _{k_1=0}^{2s}A_{k_1} + \frac{1}{U_0}\sum\limits _{l=0}^{2s}  B_{l}\prod\limits _{k_2=l+1}^{2s}A_{k_2}} \\
       & = x_0^{1 - n} \left(\frac{x_1y_0}{y_1}\right)^{n}\prod_{s = 0}^{n - 1}\frac{\prod\limits _{k_1=0}^{2s-1}C_{k_1} + x_0y_1\sum\limits _{l=0}^{2s-1}  D_{l}\prod\limits _{k_2=l+1}^{2s-1}C_{k_2}}{\prod\limits _{k_1=0}^{2s}A_{k_1} + x_1y_0\sum\limits _{l=0}^{2s}  B_{l}\prod\limits _{k_2=l+1}^{2s}A_{k_2}},
\end{align*}
\begin{align*}
x_{2n+1} & = x_1\prod_{s = 0}^{n - 1}\frac{ V_0 \prod\limits _{k_1=0}^{2s}C_{k_1} +\sum\limits _{l=0}^{2s} D_{l}\prod\limits _{k_2=l+1}^{2s}C_{k_2}}{ U_0 \prod\limits _{k_1=0}^{2s + 1}A_{k_1} +\sum\limits _{l=0}^{2s + 1}  B_{l}\prod\limits _{k_2=l+1}^{2s + 1}A_{k_2}} \\
        & = x_1\prod_{s = 0}^{n - 1}\frac{ V_0 \prod\limits _{k_1=0}^{2s}C_{k_1} +\sum\limits _{l=0}^{2s} D_{l}\prod\limits _{k_2=l+1}^{2s}C_{k_2}}{ U_0 \prod\limits _{k_1=0}^{2s + 1}A_{k_1} +\sum\limits _{l=0}^{2s + 1}  B_{l}\prod\limits _{k_2=l+1}^{2s + 1}A_{k_2}} \\
        & = x_1^{n + 1} \left(\frac{y_0}{x_0y_1} \right)^{n}\prod_{s = 0}^{n - 1}\frac{\prod\limits _{k_1=0}^{2s}C_{k_1} + x_0y_1\sum\limits _{l=0}^{2s} D_{l}\prod\limits _{k_2=l+1}^{2s}C_{k_2}}{ \prod\limits _{k_1=0}^{2s + 1}A_{k_1} + x_1y_0\sum\limits _{l=0}^{2s + 1}  B_{l}\prod\limits _{k_2=l+1}^{2s + 1}A_{k_2}},
\end{align*}
and
\begin{align*}
y_{2n} & =y_0\prod_{s = 0}^{n - 1}\frac{U_{2s}}{V_{2s +1}}\\
       & = y_0\prod_{s = 0}^{n - 1}\frac{U_0 \prod\limits _{k_1=0}^{2s-1}A_{k_1} +\sum\limits _{l=0}^{2s-1}  B_{l}\prod\limits _{k_2=l+1}^{2s-1}A_{k_2}}{ V_0 \prod\limits _{k_1=0}^{2s}C_{k_1} +\sum\limits _{l=0}^{2s} D_{l}\prod\limits _{k_2=l+1}^{2s}C_{k_2}        }\\
       & = y_0^{1 - n}\left( \frac{x_0y_1}{x_1} \right)^{n}\prod_{s = 0}^{n - 1}\frac{\prod\limits _{k_1=0}^{2s-1}A_{k_1} + x_1y_0\sum\limits _{l=0}^{2s-1}  B_{l}\prod\limits _{k_2=l+1}^{2s-1}A_{k_2}}{\prod\limits _{k_1=0}^{2s}C_{k_1} + x_0y_1\sum\limits _{l=0}^{2s} D_{l}\prod\limits _{k_2=l+1}^{2s}C_{k_2}},
\end{align*}
\begin{align*}
y_{2n+1} & =y_1\prod_{s = 0}^{n - 1}\frac{U_{2s + 1}}{V_{2s + 2}}\\
         & = y_1^{n + 1}\left( \frac{x_0}{x_1y_0} \right)^{n}\prod_{s = 0}^{n - 1}\frac{ \prod\limits _{k_1=0}^{2s}A_{k_1} + x_1y_0\sum\limits _{l=0}^{2s}  B_{l}\prod\limits _{k_2=l+1}^{2s}A_{k_2}}{\prod\limits _{k_1=0}^{2s +1}C_{k_1} + x_0y_1\sum\limits _{l=0}^{2s +1} D_{l}\prod\limits _{k_2=l+1}^{2s +1}C_{k_2} }.
\end{align*}
\section{Formulas for solutions of \eqref{1.1}}
From the previous section, the solution to \eqref{1.1} is thus given by
\begin{equation}\label{s0}
 x_{2n - 1} = x_{-1}^{1 - n} \left(\frac{x_0y_{-1}}{y_0}\right)^{n}\prod_{s = 0}^{n - 1}\frac{\prod\limits _{k_1=0}^{2s-1}c_{k_1} + x_{-1}y_0\sum\limits _{l=0}^{2s-1}  d_{l}\prod\limits _{k_2=l+1}^{2s-1}c_{k_2}}{\prod\limits _{k_1=0}^{2s}a_{k_1} + x_0y_{-1}\sum\limits _{l=0}^{2s}  b_{l}\prod\limits _{k_2=l+1}^{2s}a_{k_2}},
 \end{equation}
 \begin{equation}\label{s1}
x_{2n} = x_0^{n + 1} \left(\frac{y_{-1}}{x_{-1}y_0} \right)^{n}\prod_{s = 0}^{n - 1}\frac{\prod\limits _{k_1=0}^{2s}c_{k_1} + x_{-1}y_0\sum\limits _{l=0}^{2s} d_{l}\prod\limits _{k_2=l+1}^{2s}c_{k_2}}{ \prod\limits _{k_1=0}^{2s + 1}a_{k_1} + x_0y_{-1}\sum\limits _{l=0}^{2s + 1}  b_{l}\prod\limits _{k_2=l+1}^{2s + 1}a_{k_2}},
\end{equation}
\begin{equation}\label{s2}
y_{2n - 1} = y_{-1}^{1 - n}\left( \frac{x_{-1}y_0}{x_0} \right)^{n}\prod_{s = 0}^{n - 1}\frac{\prod\limits _{k_1=0}^{2s-1}a_{k_1} + x_0y_{-1}\sum\limits _{l=0}^{2s-1}  b_{l}\prod\limits _{k_2=l+1}^{2s-1}a_{k_2}}{\prod\limits _{k_1=0}^{2s}c_{k_1} + x_{-1}y_0\sum\limits _{l=0}^{2s} d_{l}\prod\limits _{k_2=l+1}^{2s}c_{k_2}},
\end{equation}
\begin{equation}\label{s3}
y_{2n} =  y_0^{n + 1}\left( \frac{x_{-1}}{x_0y_{-1}} \right)^{n}\prod_{s = 0}^{n - 1}\frac{ \prod\limits _{k_1=0}^{2s}a_{k_1} + x_0y_{-1}\sum\limits _{l=0}^{2s}  b_{l}\prod\limits _{k_2=l+1}^{2s}a_{k_2}}{\prod\limits _{k_1=0}^{2s +1}c_{k_1} + x_{-1}y_0\sum\limits _{l=0}^{2s +1} d_{l}\prod\limits _{k_2=l+1}^{2s +1}c_{k_2} },
\end{equation}
as long as the denominators do not vanish.\par \noindent
In the following section, we now look at the special case when all the sequences $a_n, b_n, c_n$ and  $d_n$ are constant.
\subsection{The case $a_n, b_n, c_n$ and $d_n$ are constant }
We let $a_n = a, b_n = b, c_n = c$ and $d_n = d$ where $a,b,c,d \in \mathbb{R}$. Then the solution to the system \eqref{1.1} is given by
\begin{subequations}\label{sol}
\begin{equation}\label{s0}
 x_{2n - 1} = x_{-1}^{1 - n} \left(\frac{x_0y_{-1}}{y_0}\right)^{n}\prod_{s = 0}^{n - 1}\frac{c^{2s} + dx_{-1}y_0\sum\limits _{l=0}^{2s-1}c^{l}}{a^{2s + 1} + bx_0y_{-1}\sum\limits _{l=0}^{2s} a^{l}},
 \end{equation}
 \begin{equation}\label{s1}
x_{2n} = x_0^{n + 1} \left(\frac{y_{-1}}{x_{-1}y_0} \right)^{n}\prod_{s = 0}^{n - 1}\frac{c^{2s + 1} + dx_{-1}y_0\sum\limits _{l=0}^{2s}c^{l}}{a^{2s + 2} + bx_0y_{-1}\sum\limits _{l=0}^{2s + 1} a^{l}},
\end{equation}
\begin{equation}\label{s2}
y_{2n - 1} = y_{-1}^{1 - n}\left( \frac{x_{-1}y_0}{x_0} \right)^{n}\prod_{s = 0}^{n - 1}\frac{a^{2s} + bx_0y_{-1}\sum\limits _{l=0}^{2s-1}a^{l}}{c^{2s + 1} + dx_{-1}y_0\sum\limits _{l=0}^{2s}c^{l}},
\end{equation}
\begin{equation}\label{s3}
y_{2n} =  y_0^{n + 1}\left( \frac{x_{-1}}{x_0y_{-1}} \right)^{n}\prod_{s = 0}^{n - 1}\frac{a^{2s + 1} + bx_0y_{-1}\sum\limits _{l=0}^{2s}a^{l}}{c^{2s + 2} + dx_{-1}y_0\sum\limits _{l=0}^{2s +1}c^{l}},
\end{equation}
\end{subequations}
as long as the denominators do not vanish.
\subsubsection{The case $a = c = 1$}
The solution is given by
\begin{equation}\label{s0}
 x_{2n - 1} = x_{-1}^{1 - n} \left(\frac{x_0y_{-1}}{y_0}\right)^{n}\prod_{s = 0}^{n - 1}\frac{1 + 2sdx_{-1}y_0}{1 + (2s + 1)bx_0y_{-1}},
 \end{equation}
 \begin{equation}\label{s1}
x_{2n} = x_0^{n + 1} \left(\frac{y_{-1}}{x_{-1}y_0} \right)^{n}\prod_{s = 0}^{n - 1}\frac{1 + (2s + 1)dx_{-1}y_0}{1 + (2s + 2)bx_0y_{-1}},
\end{equation}
\begin{equation}\label{s2}
y_{2n - 1} = y_{-1}^{1 - n}\left( \frac{x_{-1}y_0}{x_0} \right)^{n}\prod_{s = 0}^{n - 1}\frac{1 + 2sbx_0y_{-1}}{1 + (2s + 1)dx_{-1}y_0},
\end{equation}
\begin{equation}\label{s3}
y_{2n} =  y_0^{n + 1}\left( \frac{x_{-1}}{x_0y_{-1}} \right)^{n}\prod_{s = 0}^{n - 1}\frac{1 + (2s + 1)bx_0y_{-1}}{1 + (2s + 2)dx_{-1}y_0},
\end{equation}
where  $jbx_0y_{-1}, jdx_{-1}y_0 \neq -1$ for all $j = 1,2, \ldots, 2n$.
\subsubsection{The case $a \neq 1$ and $ c \neq 1$}
Here, \eqref{sol} simplifies to
\begin{subequations}\label{solnot1}
\begin{equation}\label{s0not1}
 x_{2n - 1} = x_{-1}^{1 - n} \left(\frac{x_0y_{-1}}{y_0}\right)^{n}\prod_{s = 0}^{n - 1}\frac{c^{2s} + dx_{-1}y_0\left( \frac{1-c^{2s}}{1-c}\right)}{a^{2s + 1} + bx_0y_{-1}\left( \frac{1-a^{2s+1}}{1-a}\right)},
 \end{equation}
 \begin{equation}\label{s1not1}
x_{2n} = x_0^{n + 1} \left(\frac{y_{-1}}{x_{-1}y_0} \right)^{n}\prod_{s = 0}^{n - 1}\frac{c^{2s + 1} + dx_{-1}y_0\left( \frac{1-c^{2s+1}}{1-c}\right)}{a^{2s + 2} + bx_0y_{-1}\left( \frac{1-a^{2s+2}}{1-a}\right)},
\end{equation}
\begin{equation}\label{s2not1}
y_{2n - 1} = y_{-1}^{1 - n}\left( \frac{x_{-1}y_0}{x_0} \right)^{n}\prod_{s = 0}^{n - 1}\frac{a^{2s} + bx_0y_{-1}\left( \frac{1-a^{2s}}{1-a}\right)}{c^{2s + 1} + dx_{-1}y_0\left( \frac{1-c^{2s+1}}{1-c}\right)},
\end{equation}
\begin{equation}\label{s3not1}
y_{2n} =  y_0^{n + 1}\left( \frac{x_{-1}}{x_0y_{-1}} \right)^{n}\prod_{s = 0}^{n - 1}\frac{a^{2s + 1} + bx_0y_{-1}\left( \frac{1-a^{2s+1}}{1-a}\right)}{c^{2s + 2} + dx_{-1}y_0\left( \frac{1-c^{2s+2}}{1-c}\right)}.
\end{equation}
\end{subequations}
\noindent \underline{\textit{The case $a = c = -1$}}\\
Then solution to the system \eqref{1.1} is given by
\begin{equation}\label{s0}
 x_{2n - 1} = x_{-1}^{1 - n} \left(\frac{x_0y_{-1}}{y_0}\right)^{n}\left(\frac{1}{-1 + bx_0y_{-1}}\right)^{n},
 \end{equation}
 \begin{equation}\label{s1}
x_{2n} = x_0^{n + 1} \left(\frac{y_{-1}}{x_{-1}y_0} \right)^{n}(-1 + dx_{-1}y_0)^{n},
\end{equation}
\begin{equation}\label{s2}
y_{2n - 1} = y_{-1}^{1 - n}\left( \frac{x_{-1}y_0}{x_0} \right)^{n}\left(\frac{1}{-1 + dx_{-1}y_0}\right)^{n},
\end{equation}
\begin{equation}\label{s3}
y_{2n} =  y_0^{n + 1}\left( \frac{x_{-1}}{x_0y_{-1}} \right)^{n}(-1 + bx_0y_{-1})^{n},
\end{equation}
where  $x_{-1}, y_0, y_{-1}, x_{0} \neq 0$ and $bx_0y_{-1}, dx_{-1}y_0 \neq 1$.

\section{Conclusion}
In this paper, we found exact solutions for a second-order system of difference equations of the form \eqref{1.1}. In the process, non-trivial symmetry generators of the system were obtained as well.
\section*{Conflict of interests}
The authors declare that there is no conflict of interests regarding the publication of this paper.

\end{document}